\documentclass{article}
 \usepackage{amsmath}
    \usepackage{amsthm}
    \usepackage{amssymb}
    \usepackage{latexsym}
\usepackage{longtable,booktabs}
\usepackage{graphicx}
\usepackage{caption}
\usepackage{mathrsfs}
\usepackage{psfragx}
\begin{document}
\newtheorem{Lemma}{Lemma}[section]
\newtheorem{Proposition}{Proposition}[section]
\newtheorem{Theorem}{Theorem}[section]
\newtheorem{Corollary}{Corollary}[section]
\newtheorem{Example}{Example}[section]
\baselineskip 0.25in
\title{\Large\bf Approximate reasoning with aggregation functions satisfying GMP rules}
\author{{Dechao Li}\thanks{Email:\ dch1831@163.com}
\qquad Qingxue Zeng\\
  {\small School of Information and Engineering},\\
{\small Zhejiang Ocean University,
  Zhoushan,
  316000, China}\\
{\small  Key Laboratory of Oceanographic Big Data Mining and Application of}\\
{\small Zhejiang Province,  Zhoushan,  316022, China}}
\date{}
\maketitle
\begin{center}
\begin{minipage}{140mm}
\begin{picture}(1,1)
 \line(1,0){400}
\end{picture}

\centerline{\bf Abstract} \vskip 3mm {\qquad
To strengthen the effectiveness of approximate reasoning in fuzzy modus ponens (FMP) and fuzzy modus tollens (FMT) problems, three approximate reasoning methods with aggregation functions are developed and their validity are investigated respectively in this paper. We firstly study some properties of fuzzy implication generated by an aggregation function. And then present an $A$-compositional rule of inference (ACRI) as an extension of Zadeh's CRI replacing $t$-norm by aggregation function. The similarity-based approximate reasoning with aggregation function (ASBR) is further discussed. Moreover, we provide the quintuple implication principle method with aggregation function (AQIP) to solve FMP and FMT problems. Finally, the validity of three approximate reasoning approaches is analyzed respectively using GMP rules in detail.}

\begin{picture}(1,1)
 \line(1,0){400}
\end{picture}
\end{minipage}
\end{center}
\vskip 4mm
\section{Introduction}
\subsection{Motivation}
\qquad Approximate reasoning has been successfully applied for model-based control, data mining, artificial intelligence,
image processing, decision making and so on. Generally speaking, approximate reasoning derives some meaningful conclusions from if-then rules and a collection of imprecise premises. Their fundamental patterns are fuzzy modus ponens (FMP) and fuzzy modus tollens (FMT) generalized from modus ponens (MP) and modus tollens (MT) in the classical logic. FMP and FMT can be represented intuitively as:
$$\textmd{Premise\ 1:}\ \textmd{IF}\ x\ \textmd{is }D\ \textmd{THEN}\ y\ \textmd{is}\ B
\qquad\qquad \textmd{Premise\ 1:}\ \textmd{IF}\ x\ \textmd{is }D\ \textmd{THEN}\ y\ \textmd{is}\ B\vspace{-2mm}$$
$$\qquad\textmd{Premise\ 2:}\ \ x\ \textmd{is }D'\qquad\quad\qquad\qquad\qquad\qquad \textmd{Premise\ 2:}\qquad\quad\qquad\qquad\ \ y\ \textmd{is }B'\qquad\vspace{-3mm}$$
\qquad\begin{picture}(1,2)
 \line(1,0){170}\qquad
 \line(1,0){170}\qquad
\end{picture}\qquad\vspace{-3mm}
$$\qquad\quad\ \textmd{Conclusion:}\qquad\quad\qquad\qquad y\ \textmd{is}\ B',\qquad\qquad \textmd{Conclusion:}\ x\ \textmd{is}\ D',\qquad\qquad\qquad\qquad\qquad$$
where $D$ and $D'$ are fuzzy sets on the universe $U$ while $B$ and $B'$ are fuzzy sets on the universe $V$.

To obtain $B'(D')$ from $B(D)$, the compositional rule of inference (CRI) method was proposed by Zadeh in 1973 \cite{Zadeh}. In Zadeh's CRI, Premise 1 is translated into a fuzzy relation $R$ using Zadeh implication. Then $B'(D')$ is calculated
by combining $D'(B')$ and fuzzy relation $R$ with the sup-min composition.
After, the general CRI methods for FMP and FMT are developed as follows:
$$B'(y)=\bigvee_{x\in U}D'(x)\ast(D(x)\rightarrow B(y)),$$
$$D'(x)=\bigvee_{y\in V}B'(y)\ast(D(x)\rightarrow B(y)),$$
where $\ast$ is a t-norm, $\rightarrow$ is a fuzzy implication. Instead of t-norm, Ruan and Kerre also extended the CRI method by $n$-ary operator $T_n$ \cite{Ruan}. Moreover, Cappelle et al. studied the CRI method in case where a binary function $F$ on [0,1] is used to explain Premise 1 (that is, Premise  1 is translated into $F(D(x), B(y))$) \cite{Cappelle}. After, Koles\'{a}rov\'{a} and Kerre investigated the CRI method in special case where the function $F$ is a t-norm \cite{Kolesarova}.

It is necessary to mention that Trillas et al. represented the following FMT \cite{Trillas}:
$$\textmd{(Trillas'\ FMT)}\qquad\qquad\textmd{Premise\ 1:}\ \textmd{IF}\ x\ \textmd{is }D\ \textmd{THEN}\ y\ \textmd{is}\ B\qquad\qquad\qquad\qquad\qquad\ \vspace{-2mm}$$
$$\qquad\textmd{Premise\ 2:}\quad\qquad\qquad\quad  y\ \textmd{is\ not }B'\qquad\vspace{-3mm}$$
\qquad\begin{picture}(1,2)
 \quad\qquad\qquad\qquad\qquad\quad\line(1,0){170}\qquad
\end{picture}\qquad\vspace{-1mm}
$$\qquad\quad\ \textmd{Conclusion:}\ x\ \textmd{is\ not}\ D'.\quad\qquad\qquad\qquad\quad$$
To solve Trillas' FMT, the strong negation is used to explain the connective ``not". With CRI method, the conclusion is obtained by Trillas et al. as follows $$N(D'(x))=\bigvee_{y\in V}N(B'(y))\ast(D(x)\rightarrow B(y)).$$
We do not consider Trillas' FMT in the rest of this work.

Although CRI method is simple in computation, there are still some deficiencies in CRI method as pointed out by
some researchers \cite{Baldwin,Mizumoto,Turksen,Wang,Zhou}. To overcome these deficiencies, Turksen and Zhong suggested similarity-based approximate reasoning (SBR) method which does not require to construct the fuzzy relation \cite{Turksen}. After, Raha et al. developed an SBR method using a new measure for similarity between two fuzzy sets \cite{Raha}.
In order to provide a logical foundation for FMP and FMT problems, Wang and Pei proposed triple implication principle (TIP) for fuzzy reasoning \cite{Pei,Wang1}. To improve the quality of TIP method, Zhou et al. investigated quintuple implication principle (QIP) for FMP and FMT problems \cite{Zhou}. Most importantly, it is found that Mamdani-type fuzzy inference is same as fuzzy inference with QIP method using G\"{o}del implication.

To measure the validity of inference scheme to solve the FMP and FMT problems, Magrez and Smets \cite{Magrez} proposed some commonly accepted axioms (Also inferred as GMP rules) in the following:

(GMP1)\quad $B\subseteq B'$;

(GMP2)\quad If $D'\subseteq D''$, then $B'\subseteq B''$;

(GMP3)\quad If $D'=D^C$, then $B'=V$, where $D^C$ is the complement of $D$;

(GMP4)\quad If $D'=D$, then $B'=B$.

In order to make better use of approximate reasoning, it becomes a core topic to measure the validity of inference scheme using GMP rules \cite{Baets, Cornelis, Cornelis1, Mas, Mas2}.

It is well known that the results of approximate reasoning depend completely
on the choice of logical connectives. However, as some researchers \cite{Bustince,Fodor} pointed out the associativity or commutativity of the connectives ``and" and ``or" is not demanded in classification problems and
decision making. Considered aggregation functions play an important role in decision making and fuzzy logic, aggregation functions are a better substitute for the t-norms and t-conorms by in the actual classification problems and decision making.

Moreover, fuzzy
implication, as an important logical connective, is
used to formalize ``if ... then" rule in fuzzy system. There
exist many families of fuzzy implications, such as well-known R-, S- and QL-implications, $f$- and $g$-implications, probabilistic implication, probabilistic S-implication and so on. According to the generation methods of fuzzy implications, fuzzy implications can be classified into two types as follows: i. generated by the binary functions on [0, 1], such as R-, (S, N)-, QL-implications, residual implications derived from overlap functions and probabilistic implications \cite{Baczynski,Dimuro,Dimuro1,Grzegorzewski}; ii. generated by the unary functions on [0, 1], for instance, $f$- and $g$-implications \cite{Yager}. As the t-norm and t-conorm are two special aggregation functions (See Definition 2.5), it is very interesting topic to investigate fuzzy implications generated by aggregation functions. As mentioned above, the actual classification problems and decision making also trigger us to study the fuzzy implications generated by aggregation functions. Thus, our motivation is to develop three approximate reasoning approaches using aggregation functions and fuzzy implications generated by them. Most importantly, we will pay close attention to the validity of three approximate reasoning methods.
\subsection{Contribution of this research}
\qquad As is well known that the FMP and FMT are two models to obtain the conclusion from imprecise premises. They also play a pivotal role in decision making. Therefore, it is not difficult to see that more options of fuzzy implications and aggregation functions result in more flexibility in decision making. Based on the discussion above, we mainly develop three approximate reasoning approaches with aggregation functions to solve FMP and FMT problems in this paper. And what's more, the validity of three approximate reasoning methods is respectively discussed using GMP rules.  We first investigate some properties of fuzzy implication generated by an aggregation function. Based on such fuzzy implication and aggregation function, three approximate reasoning approaches are developed to solve FMP and FMT problems.
In a word, the contributions of this paper include:

(1) To study the properties of fuzzy implication generated by an aggregation function.

(2) To construct three approximate reasoning methods using aggregation functions (that is, ACRI method, ASBR method and AQIP method).

(3) To investigate the validity of these three approximate reasoning methods using GMP rules.

This paper is composed as follows. In Section 2,
some definitions of basic notions and notations are presented. Section 3 studies some properties of fuzzy implication generated by an aggregation function. In Section 4, the ACRI method with aggregation function is discussed. In Section 5,  we propose the ASBR method. Section 6 provides
the AQIP method for FMP and FMT problems.
\section{Preliminary}
\qquad In order to make this work more self-contained, we introduce the main concepts and properties
employed in the rest of the paper.
\subsection{Negation, aggregation function and fuzzy implication}
{\bf Definition 2.1}\cite{Lowen} A function $N: [0,1]\rightarrow [0,1]$ is called
a fuzzy negation if

(N1) $N(0)=1, \ N(1)=0$;

(N2) $N(x)\geq N(y)\ \textmd{if}\ x\leq y,\ \forall\ x, y\in [0,1]$.

 Further, a fuzzy negation $N$ is strict if it satisfies the following properties:

(N3) $N$ is continuous;

(N4) $N(x)> N(y)\ \textmd{if}\ x<y$.

A fuzzy negation is strong if it is involutive, i.e.,

(N5)  $N(N(x))=x,
\forall\ x\in [0,1]$.
\\{\bf Example 2.2}\cite{Lowen} The negation $N_0(x)=1-x$ is strong. It also is called the standard negation.
\\{\bf Definition 2.3}\cite{Grabisch} A function $A:[0,1]^n\rightarrow [0,1]$ is said to be an $n$-ary aggregation function if the following statements hold:

(A1) $A$ satisfies the boundary conditions: $A(0, 0,\cdots, 0)=0$ and $A(1, 1,\cdots, 1)= 1$;

(A2) $A$ is non-decreasing in each variable.
\\{\bf Definition 2.4}\cite{Grabisch} Let $A$ be a binary aggregation function.

i. An element $a\in[0,1]$ is said to be a left (right) annihilator if $A(a,x)=a\ (A(x,a)=a)$ for any $x\in[0,1]$; $a$ is an annihilator if $A(a,x)=A(x,a)=a$ for any $x\in[0,1]$;

ii. $e\in[0,1]$ is said to be a left (right) neutral element if  $A(e, x)=x\ (A(x, e)=x)$ for any $x\in[0,1]$; $e\in[0,1]$ is a neutral element if  $A(e, x)=A(x, e)=x$ for any $x\in[0,1]$.
\\{\bf Definition 2.5}\cite{Grabisch} A binary aggregation function $A$ is said to be

i. symmetric or commutative if $A(x,y)=A(y,x)$ for any $x,y\in[0,1]$;

ii. associative if $A(x,A(y,z))=A(A(x,y),z)$ for any $x,y,z\in[0,1]$;

iii. conjunctive if $A\leq \textmd{min}$;

iv. disjunctive if $A\geq \textmd{max}$;

v. averaging if $\textmd{min}\leq A\leq \textmd{max}$;

vi. a semi-copula if 1 is a neutral element;

vii. a t-norm if it is an associative and commutative semi-copula;

viii. dual to a semi-copula if 0 is a neutral element;

ix. a t-conorm if it is dual to a t-norm;

x. a copula if it is a semi-copula which is two-increasing, i.e., $A(x_1,y_1)-A(x_1,y_2)-A(x_2,y_1)+A(x_2,y_2)\geq0$ holds for all $x_1,y_1,x_2,y_2\in [0,1]$ such that $x_1\leq x_2$ and $y_1\leq y_2$.
\\{\bf Definition 2.6}\cite{Baczynski}  A fuzzy implication is a function  $I: [0, 1]^2\rightarrow[0, 1]$ which satisfies for any $x,y,z\in [0,1]$:

(I1) Non-increasing in the first variable, i.e., if $x\leq y$ then $I(x,z)\geq I(y, z)$;

(I2) Non-decreasing in the second variable, i.e., if $y\leq z$ then $I(x,y)\leq I(x, z)$;

(I3) $I(0,0)=1$;

(I4) $I(1,1)=1$;

(I5) $I(1,0)=0$.

By Definition 2.6, we directly obtain the fact that a fuzzy implication satisfies the following properties:

(LB) Left boundary condition, $I(0, y)= 1,\forall\ y\in[0, 1]$;

(RB) Right boundary condition, $I(x, 1)=1,\forall\ x\in[0, 1]$.
\\{\bf Definition 2.7}\cite{Baczynski} A fuzzy implication  $I : [0, 1]^2\rightarrow[0, 1]$ satisfies:

(NP) Left neutrality property, if $I(1, y)=y,\forall\ y\in[0, 1]$;

(IP) Identity principle, if $I(x, x) = 1, \forall\ x\in[0, 1]$;

(EP) Exchange principle, if $I(x, I(y, z)) = I(y,
I(x, z)), \forall\ x, y, z\in[0, 1]$;

(CP(N)) Law of contraposition with a fuzzy negation $N$ if $I(x, y) = I(N(y),N(x)), \forall\ x, y\in[0, 1]$;

(OP) Ordering property, if $I(x, y) =
1\Longleftrightarrow x\leq y, \forall\ x, y\in[0, 1]$.

Different classes of implications can be found in many literatures. Among
them we only emphasize the following classes of fuzzy implications.
 \\{\bf Definition 2.8}\cite{Baczynski} An R-implication is a function $I_T:[0,1]^2\rightarrow[0,1]$ associated with a t-norm $T$ defined by $I_T(x,y)=\sup\{z|T(x,z)\leq y\}$.
  \\{\bf Definition 2.9}\cite{Pradera}  An $(A,N)$-implication is a function $I_{A,N}:[0,1]^2\rightarrow[0,1]$ associated with a disjunctor (that is, an aggregation function having an annihilator 0) $A$ and a fuzzy negation $N$ defined by $I_{A,N}(x,y)=A(N(x), y)$.
\\{\bf Definition 2.10}\cite{Yager} Let $f:[0,1]\rightarrow [0,\infty]$ be a strict decreasing and continuous mapping with
$f(1)=0$. An $f$-generated implication, which is a function $I_f:[0,1]^2\rightarrow[0,1]$ with an $f$-generator, is defined by
$I_{f}(x,y)=f^{(-1)}(xf(y))$ with the understanding that $0\times \infty=0$,\vspace{2mm} where $f^{(-1)}$ is pseudoinverse of $f$ defined as $f^{(-1)}(x)=\left\{\begin{array}{ll}
                                                                                    f^{-1}(x) & x\leq f(0) \\
                                                                                    0 & \textmd{otherwise}
                                                                                  \end{array}\right.$.\vspace{2mm}
\\{\bf Definition 2.11}\cite{Yager} Let $g:[0,1]\rightarrow [0,\infty]$ be a strict increasing and continuous mapping with
$g(0)=0$. A $g$-generated implication, which is a function $I_g:[0,1]^2\rightarrow[0,1]$ with a $g$-generator, is defined by
$I_{g}(x,y)=g^{(-1)}\left(\frac{g(y)}{x}\right)$ with the understanding that $0\times \infty=\infty$,  where $g^{(-1)}$ is pseudoinverse of $g$.\vspace{1mm}
\\{\bf Definition 2.12}\cite{Grzegorzewski} Let $C$ be a copula. A function $I_C:[0, 1]^2\rightarrow[0, 1]$  given by $I_C(x,y)\vspace{1mm}=\left\{\begin{array}{ll}
                                                                                   \frac{C(x,y)}{x} & x>0 \\
                                                                                    1 & \textmd{otherwise}
                                                                                  \end{array}\right.$
is called a probabilistic implication (based on a copula $C$).\vspace{1mm}
\\{\bf Definition 2.13}\cite{Grzegorzewski} Let $C$ be a copula. A function $\tilde{I}_C:[0, 1]^2\rightarrow[0, 1]$ given by
$\tilde{I}_C(x,y)=C(x,y)-x+1$
is called a probabilistic S-implication (based on a copula $C$).
\subsection{Raha's similarity-based approximate reasoning}
\qquad Similarity-based approximate reasoning methods can be founded in many literatures. In this subsection, we only recall the similarity-based approximate reasoning method proposed by Raha et al. in \cite{Raha}.
Let $F(U)$ denote all fuzzy sets defined on the universe $U$.
\\{\bf Definition 2.14}\cite{Raha} A function $S:F(U)\times F(U)\rightarrow [0, 1]$ is called a similarity measure if it satisfies the following properties for any $D, D'\in F(U)$:

(S1) $S(D, D')=S(D', D)$;

(S2) $S(D, D')=1$ if and only if $D=D'$;

(S3) $D, D'$ are simultaneously not null, that is,  $\min(D(x), D'(x))=0$ for all $x\in U$ if $S(D,D')=0$;

(S4) $S(D, D'')\leq \min(S(D, D'),S(D',D''))$ if $D\subseteq D'\subseteq D''$.

In order to obtain the conclusion in FMP problem,  Raha et al. presented a novel similarity-based approximate reasoning method  \cite{Raha}. In their proposed method, Premise 1 is interpreted as a conditional fuzzy relation $R(D, B)$ while the conclusion is interpreted as a modified conditional relation $R(D, B|D')$. And an algorithm for similarity-based approximate reasoning is shown as follows:

Step 1. Translate premise 1 and compute $R(D, B)$ using some suitable translating rules (possibly, a t-norm
operator).

Step 2. Compute $S(D', D)$ between the fact $D'$ and the antecedent $D$ using some suitable similarity measures.

Step 3. Modify $R(D, B)$ with $S(D', D)$ to obtain the modified conditional relation $R(D, B|D')$ using some schemes.

Step 4. Use the sup-projection operation on $R(D, B|D')$ to obtain $B'$ as $$B'(y)=\mathop{\sup}\limits_{x\in U}R(D,B|D')(x,y).$$

In order to compute the conclusion $R(D, B|D')$, the following axioms are proposed:

(AX1) If $S(D', D)=1$, then $R(D, B|D')(x, y)=R(D, B)(x, y)$;

(AX2) If $S(D', D)=0$, then $R(D, B|D')(x, y)=1$;

(AX3) As $S(D', D)$ increases from 0 to 1, $R(D, B|D')(x, y)$ decreases uniformly from 1 to $R(D, B)(x, y)$, that is,  $R(D, B|D')\supseteq R(D, B)$ holds for any $D'\in F(U)$.

Then $R(D, B)$ is constructed in following ways:

Case 1. When $R(D, B)(x, y)=T(A(x), B(y))$, where $T$ is a t-norm.

Case 2. When $R(D, B)(x, y) =I_T(A(x), B(y))$, where $I_T$ is an R-implication.

Finally, Raha et al. obtained the conclusions $B'_1$ and $B'_2$  as
$$B'_1(y)=\mathop{\sup}\limits_{x\in U}I_T(S(D,D'),T(D(x),B(y))),$$
$$B'_2(y)=\mathop{\sup}\limits_{x\in U}I_T(S(D,D'),I_T(D(x),B(y))).$$
\subsection{Quintuple implication principle for FMP and FMT}
\qquad Zhou et al. proposed the following quintuple implication principle (QIP) for solving FMP and FMT problems \cite{Zhou}. In \cite{Li1}, Li and Qin extended the quintuple implication principle for FMP and FMT as follows.
\\{\bf Quintuple implication principle for FMP} Let $D,D'\in F(U)$ and  $B\in F(V)$. Suppose the maximum of following formula
\begin{equation}M(x,y)=I(I(D(x),B(y)), I(I(D'(x),D(x)),
I(D'(x),B'(y))))
\end{equation}
exists for every $x\in U$ and $y\in V$, where $I$ is a fuzzy implication on [0,1]. The solution $B'$ of FMP should be the smallest fuzzy subset on $V$ such that Eq.(1) takes its maximum.
\\{\bf Quintuple implication principle for FMT}  Let $D\in F(U)$ and  $B, B'\in F(V)$. Suppose the maximum of following formula
\begin{equation}N(x,y)=I(I(D(x),B(y)), I(I(B(y),B'(y)),
I(D(x),D'(x))))
\end{equation}
exists for every $x\in U$ and $y\in V$. The solution $D'$ of FMT is the smallest fuzzy subset on $U$ such that Eq.(2) takes its maximum.
\\{\bf Lemma 2.15}\cite{Li1} i. If $I$ satisfies (I2), then the greatest value of the formulas (1) and (2) are
$$\max_{x\in U,y\in V}M(x,y)=I(I(D(x),B(y)), I(I(D'(x),D(x)),
I(D'(x),1)))$$\vspace{-2mm} and
 $$\max_{x\in U,y\in V}N(x,y)=I(I(D(x),B(y)), I(I(B(y),B'(y)),
I(D(x),1))).$$

ii. Moreover, if $I$ is right-continuous with respect to the second variable, then the QIP solution of FMP (FMT) exists and is unique.
\\{\bf Theorem 2.16}\cite{Zhou} Suppose $I$ is an R-implication induced by a left-continuous t-norm $T$. Then the QIP solutions of FMP and FMT are as follows:
$$B'(y)=\sup_{x\in U}T(D'(x),T(I(D'(x),D(x)),I(D(x),B(y)))),$$
$$D'(x)=\sup_{y\in V}T(D(x),T(I(D(x),B(y)),I(B(y),B'(y)))).$$
\section{Residual implication generated by an aggregation function}
\qquad In this section, we will investigate some properties of fuzzy implication generated by an aggregation function.

 Let $A$ be a function from $[0,1]^2$ to $[0,1]$.  Then, we can define a function $I_A :[0,1]^2\rightarrow [0,1]$ as
\begin{equation}I_A(x,y)=\sup\{z\in [0,1]|A(x, z)\leq y\}, \forall x,y \in [0,1].\end{equation}

For an aggregation function $A$, we have the following result.
\\{\bf Lemma 3.1} Let $A$ be an aggregation function. Then the following statements are equivalent:

i. $A$ is left-continuous with respect to the second variable;

ii. $A$ and $I_A$ defined in Eq.(3) satisfy the residuation property (RP), i.e.
$$A(x,z)\leq y\Longleftrightarrow z\leq I_A(x,y), \ \forall x,y,z\in[0,1];\eqno(\textmd{RP})$$

iii. $I_A(x,y)=\max\{z\in [0,1]|A(x, z)\leq y\}, \forall x,y \in [0,1]$.
\\{\bf Proof.} This proof is similar to that of Proposition 2.5.2 in \cite{Baczynski} and Theorem 2 in \cite{Krol}.
\\{\bf Remark 1.}  The above proposition also appeared in \cite{Demirli, Jayaram} when $A$ is a semi-copula. In \cite{Krol}, Kr\'{o}l considered the case where $A$ is a conjunctor. However, it is sufficient to demand that $A$ is an aggregation function here.

Considering an aggregation function $A: [0,1]^2\rightarrow [0,1]$ under certain conditions, it is possible to define a class of
fuzzy implications according to Eq.(3). In \cite{Ouyang}, it is proved that $I_A$ is a fuzzy implication if the aggregation function $A$ satisfies the following conditions:
\begin{equation} A(1,y)>0\ \textmd{for\ any}\ y>0, \end{equation}
\begin{equation} A(0,y)=0\ \textmd{for\ any}\ y<1.\end{equation}
In this case, we say $I_A$ is a residual implication induced by the aggregation function $A$ (for short, R-implication). Notice that the above result also appeared in \cite{Krol}.

We try to obtain an aggregation function from a fuzzy implication in turn. Let $I$ be a fuzzy implication. The function $A_I :[0,1]^2\rightarrow [0, 1]$ is defined by:
\begin{equation}
A_I(x,y)=\inf\{z \in[0, 1]|I(x,z)\geq y\}, \forall x,y\in[0,1].\end{equation}
\qquad Similar to Lemma 3.1, we have the following result.
\\{\bf Lemma 3.2}\cite{Baczynski,Krol} Let $I$ be a fuzzy implication. Then, the following statements are equivalent:

i. $I$ is right-continuous with respect to the second variable;

ii. $A_I$ defined Eq.(6) and $I$ satisfy the residuation property, i.e.
 $$A_I(x,z)\leq y\Longleftrightarrow z\leq I(x,y),\  \forall x,y,z\in[0,1];\eqno(\textmd{RP*})$$
\quad\ iii. $A_I(x,y)=\min\{z\in [0,1]|I(x, z)\geq y\}, \forall x,y \in [0,1]$.
\\{\bf Lemma 3.3} If the fuzzy implication $I$ satisfies the condition $I(1,y)<1$  for all $y\in[0,1)$, then the function $A_I$ defined in Eq.(6) is an aggregation function.
\\{\bf Proof.} Obviously, $0\in \{z|I(0,z)\geq 0\}$ holds. This implies $A_I(0,0)=0$. Similarly, we have $A_I(1,1)=\inf\{z|I(1,z)=1\}=1$ by the condition $I(1,y)<1$  for all $y\in[0,1)$.

 It is not difficult to obtain the fact that $A_I$ is nondecreasing in two variables since a fuzzy implication $I$ satisfies (I1) and (I2).
 \\{\bf Remark 2.} Indeed, the above lemma also appeared in \cite{Krol}. In this case, 0 is an annihilator of the aggregation function $A_I$.

 As an extension of the Theorem 2.5.14 in \cite{Baczynski}, we further get the following statement.
 \\{\bf Theorem 3.4} Let fuzzy implication $I$  be right-continuous with respect to the second variable. Then $I=I_{A_I}$, i.e., $I(x,y)=\max\{z|A_I(x,z)\leq y\}$ for any $x,y\in[0,1]$, where the function $A_I$ is defined in Eq.(6).
\\{\bf Proof.} We firstly verify that $I_{A_I}$ is a fuzzy implication. Obviously, $I_{A_I}$ satisfies (I1) and (I2) by the definition of $A_I$. Therefore, it is sufficient to verify $I_{A_I}(0,0)=I_{A_I}(1,1)=1$ and $I_{A_I}(1,0)=0$. Since $I$ is right-continuous with respect to the second variable, $A_I(1,z)\leq 1\Longleftrightarrow z\leq I(1,1)=1$ holds for all $z\in[0,1]$ by Lemma 3.2.  This implies that   $I_{A_I}(1,1)=\sup\{z|A_I(1,z)\leq1\}=1$.

Since $I$ satisfies (RB), $I(0,y)\geq z$ holds for all $y,z\in[0,1]$. According to the definition of $A_I$, we have $I_{A_I}(0,0)=\sup\{0|A_I(0,z)=0\}=1$.

Assume that $A_I(1,z)=0$. That is, $\min\{y|I(1,y)\geq z\}=0$. The right-continuity of $I$ with respect to the second variable implies $z=0$. This implies $I_{A_I}(1,0)=0$.

Next, we prove $I=I_{A_I}$. Since $I(x,y)\leq I(x,y)$ and $A_I(x,I(x,y))\leq y$ hold for all $x,y\in[0,1]$, $I(x,y)\leq I_{A_I}(x,y)$ holds.

On the other hand, we can assert that $A_I$ is left-continuous with respect to the second variable. Indeed, let any $x, y_i\in [0,1]$ and $i\in S$. $A_I(x, \mathop{\vee}\limits_iy_i)\geq \mathop{\vee}\limits_iA_I(x,y_i)$ holds for every $i\in S$. Let $\mathop{\vee}\limits_iA_I(x,y_i)=y$. Then we have $A_I(x, y_i)\leq y$ for every $i\in S$. According to Lemma 3.2, $y_i\leq I(x,y)$ holds for every $i\in S$. This implies that $\mathop{\vee}\limits_iy_i\leq I(x,y)$. Again, we obtain
$A_I(x, \mathop{\vee}\limits_iy_i)\leq y$ by Lemma 3.2.

Obviously, the left-continuity of $A_I$ with respect to the second variable implies that
$I_{A_I}(x,y)\geq I_{A_I}(x,y)\Longleftrightarrow A_I(x,I_{A_I}(x,y))\leq y$ holds for any $x,y\in [0,1]$. Since $A_I(x,z)\leq A_I(x,z)$, we have $I(x, A_I(x,z))\geq z$ for any $x, z\in [0,1]$ by (RP*). Especially, take $z=I_{A_I}(x,y)$. Then we obtain $I(x, A_I(x,I_{A_I}(x,y)))\geq I_{A_I}(x,y)$. This implies $I(x,y)\geq I(x, A_I(x,I_{A_I}(x,y)))\geq I_{A_I}(x,y)$.

By the discussion above, we get $I(x,y)=I_{A_I}(x,y)$.
\\{\bf Remark 3.} i. In \cite{Krol}, the fuzzy implication $I$ not only is right-continuous with respect to the second variable but also satisfies
the condition $I(1, y)<1$ for all $y\in[0,1)$.

ii. Theorem 2.5.14 in \cite{Baczynski} demands that $I$ satisfies (I2), (EP), (OP) and is right-continuous with respect to the second variable. In this case, $A_I$ is a t-norm.

iii. The above result shows that all right-continuous with respect to the second variable fuzzy implications (including well-known R-, S- and QL-implications, $f$- and $g$-implications, probabilistic implications, probabilistic S-implications, etc.) can be obtained as R-implications induced by aggregation functions.
\section{$A$-compositional rule of inference satisfying GMP rules}
\subsection{$A$-compositional rule of inference with aggregation function}
\qquad In this subsection, we study the composition rule of inference method based on
the aggregation function $A$ satisfying GMP rules.
\\{\bf Definition 4.1} Let $R$ and $S$ be two fuzzy relations on $U\times V$ and $V\times W$, respectively. A $\sup-A$ composition
of the fuzzy relations $S$ and $R$ is defined as a relation $S\circ_{A} R$ on $U\times W$ in the following:
\begin{equation}
(S\circ_A R)(x,z)=\sup_{y\in V}A(S(x,y),R(y,z)).
\end{equation}

Based on it, the ACRI methods for FMP and FMT problems can be developed as follows:
\begin{equation}
B'(y)=\bigvee_{x\in U}A(D'(x),I(D(x), B(y))),\end{equation}
\begin{equation}D'(x)=\bigvee_{y\in V}A(B'(y),I(D(x), B(y))),
\end{equation}
where $A$ is an aggregation function and $I$ a fuzzy implication.

Next, we shall look for the aggregation functions which the ACRI methods for FMP and FMT problems satisfy GMP rules for an arbitrary fixed fuzzy implication $I$.
\\{\bf Theorem 4.2} Let $I$ be a fuzzy implication and $f(y)=I(1,y)$ a strictly increasing function on $[0,1]$. Then there exists an aggregation function $A_I$ defined as Eq.(6) such that the following statements of ACRI hold:

i. ACRI satisfies (GMP1) if $D'$ is normal;

ii. ACRI satisfies (GMP2);

iii. ACRI satisfies (GMP3) if $D^C$ is normal;

iv. ACRI satisfies (GMP4) if $D$ is normal.
\\{\bf Proof.} We only consider the ACRI method for FMP. The ACRI for FMT can be considered similarly.
 i. Since $D'$ is normal, there exists $x_0\in U$ such that $D'(x_0)=1$. And then $B'(y)=\mathop{\bigvee}\limits_{x\in U}A_I(D'(x),I(D(x), B(y)))\geq A_I(D'(x_0),I(D(x_0), B(y)))\geq A_I(1,I(1,B(y)))=$\vspace{1mm} $\inf\{z\in[0,1]|
I(1,z)\geq I(1,B(y))\}=B(y)$ holds, where we use the strict increase of $f(y)=I(1,y)$.

 ii. By Lemma 3.3, we can immediately get the fact that ACRI satisfies (GMP2).

iii. Since $D^C$ is normal, there exists $x_0$ such that $D^C(x_0)=1$. This implies $B'(y)=\mathop{\bigvee}\limits_{x\in U}A_I(D^C(x),I(D(x), B(y)))\geq A_I(D^C(x_0),I(D(x_0), B(y)))=A_I(1,I(0, B(y)))=$\vspace{1.5mm} $A_I(1,1)=1$.

iv. Let $D'=D$. This means $B'(y)=\mathop{\bigvee}\limits_{x\in U}A_I(D(x),I(D(x), B(y)))$. Since $D$ is normal,\vspace{1mm} we have $B(y)=A_I(1,I(1, B(y)))\leq \mathop{\bigvee}\limits_{x\in U}A_I(D(x),I(D(x), B(y)))\leq B(y)$. Thus, $B=B'$.\vspace{1mm}
\\{\bf Corollary 4.3}  Let $I$ be an $(A,N)$-implication and $f(y)=A(0,y)$ a strictly increasing function on $[0,1]$. Then there exists an aggregation function $A_I$ defined as Eq.(6) such that the following statements of ACRI hold:

i. ACRI satisfies (GMP1) if $D'$ is normal;

ii. ACRI satisfies (GMP2);

iii. ACRI satisfies (GMP3) if $D^C$ is normal;

iv. ACRI satisfies (GMP4) if $D$ is normal.
\\{\bf Proof.} For an $(A,N)$-implication, we have $f(y)=I(1,y)=A(0,y)$. Then the results can be proved similarly.
\\{\bf Remark 4.} As showed in \cite{Pradera}, all fuzzy implications can be obtained as $(A,N)$-implications. This means that we can construct
ACRI method using the aggregation function $A$ according to Corollary 4.3.

Similarly, we can verify the following results.
\\{\bf Corollary 4.4}  Let $I$ be an $f$-implication or $g$-implication. Then there exists an aggregation function $A_I$ defined as Eq.(6) such that the following statements of ACRI hold:

i. ACRI satisfies (GMP1) if $D'$ is normal;

ii. ACRI satisfies (GMP2);

iii. ACRI satisfies (GMP3) if $D^C$ is normal;

iv. ACRI satisfies (GMP4) if $D$ is normal.
\\{\bf Corollary 4.5}  Let $I$ be a probabilistic implication or probabilistic S-implication and $f(y)=C(1,y)$ a strictly increasing function on $[0,1]$. Then  there exists an aggregation function $A_I$ defined as Eq.(6) such that the following statements of ACRI hold:

i. ACRI satisfies (GMP1) if $D'$ is normal;

ii. ACRI satisfies (GMP2);

iii. ACRI satisfies (GMP3) if $D^C$ is normal;

iv. ACRI satisfies (GMP4)  if $D$ is normal.

In turn, we look for the fuzzy implications which the ACRI methods for FMP and FMT problems satisfy GMP rules for an arbitrary fixed aggregation function $A$.
\\{\bf Theorem 4.6} Let $A$ be a left-continuous with respect to the second variable aggregation function. If $A$ has a left neutral element 1 and satisfies Eq.(5). Then there exists a fuzzy  implication $I_A$ defined as Eq.(3) such that the following statements of ACRI method based on $A$ hold:

i. ACRI satisfies (GMP1) if $D'$ is normal;

ii. ACRI satisfies (GMP2);

iii. ACRI satisfies (GMP3) if $D^C$ is normal;

iv. ACRI satisfies (GMP4)  if $D$ is normal.
\\{\bf Proof.} This result comes from Lemma 3.1.
\\{\bf Theorem 4.7}  Let $A$ be a left-continuous with respect to the second variable aggregation function and $I$ a fuzzy implication. If $A$ has a left neutral element 1 and $I$ satisfies (NP), then the ACRI based on $A$ and $I$ satisfies (GMP1)-(GMP4) if and only if $I\leq I_A$.
\\{\bf Proof.} $(\Longrightarrow)$ Since the ACRI method satisfies GMP4, $A(x,I(x,y))\leq y$ holds for any $x,y\in[0,1]$. By Lemma 3.2, we obtain $I(x,y)\leq I_A(x,y)$.

$(\Longleftarrow)$ Let us verify the ACRI method satisfies (GMP1)-(GMP4).

 i. Since $D'$ is normal, there exists $x_0\in U$ such that $D'(x_0)=1$. And then $B'(y)=\mathop{\bigvee}\limits_{x\in U}A(D'(x),I(D(x), B(y)))\geq A(D'(x_0),I(D(x_0), B(y)))\geq A_I(1,I(1,B(y))=B(y)$ holds.\vspace{1mm}

 ii. Obviously, ACRI method satisfies (GMP2).

iii. Since $D^C$ is normal, there exists $x_0$ such that $D^C(x_0)=1$. This implies $B'(y)=\mathop{\bigvee}\limits_{x\in U}A(D^C(x),I(D(x), B(y)))\geq A(D^C(x_0),I(D(x_0), B(y)))=A(1,I(0, B(y)))\geq A_I(1,1)=1$.

iv. Let $D'=D$. In this case $B'(y)=\mathop{\bigvee}\limits_{x\in U}A(D(x),I(D(x), B(y)))$. Since $D$ is normal, we\vspace{1mm} have $B(y)=A(1,I(1, B(y)))\leq \mathop{\bigvee}\limits_{x\in U}A(D(x),I(D(x), B(y)))\leq \mathop{\bigvee}\limits_{x\in U}A(D(x),I_A(D(x), B(y)))\leq B(y)$. Thus, $B'=B$.
\subsection{Approximate reasoning in ACRI method with multiple fuzzy rules}
\qquad We have studied ACRI method for a single fuzzy rule above. In practical applications, it needs to deal with approximate reasoning with multiple fuzzy rules. Therefore, this subsection extends the ACRI method in the case of multiple fuzzy rules involved.
It is well known that IF-THEN rule base is the main parts of fuzzy system. And the fuzzy rule base of multiple-input and single-output (MISO) fuzzy system consists of rules as follows:
\begin{equation}\textmd{R}_j:\ \ \ \textmd{IF} \ \ x_{1}\ \textmd{is}\ D^1_j\ \textmd{ AND} \ \ x_{2}\ \textmd{is}\ D^2_j\ \textmd{ AND}
 \cdots\ \textmd{AND} \ x_{m}\ \textmd{is}\ D^m_j\ \textmd{THEN  }\ y\ \textmd{is} \ B_{j}, \qquad
\end{equation}
where $x_i(i=1,2,\cdots, m)$ and $y$ are variables and $D^i_{j}(j=1,2,\cdots, n)$ and $B_j$ are specific linguistic expressions  expressing properties of values of
 $x_i$ and $y$, respectively.

 Let $D_{j}=D^{1}_{j}\times D^{2}_{j}\times \cdots\times
D^{m}_{j}$ and $\mathbf{x} = (x_1,x_2,\cdots,x_m)$.  Then each fuzzy rule R$_j$ can be regarded as a fuzzy relation $R_j$ with a membership function
$R_j(\mathbf{x},y)=D_j(\mathbf{x})\rightarrow B_j(y)$. Further,
 t-norms are employed to
  evaluate the ANDs in the fuzzy rules. In order to obtain the result of inference $B'$ from an input and fuzzy rules, we employ two schemes in general \cite{Wang}. One is First Infer Then Aggregate (FITA). That is, for a given input $D'$, we first compose $D'$ with each fuzzy rule to infer $m$ individual $B_j'$. And then aggregate $B_j'$ into the overall output $B'_{FITA}$. In this case, the output can be written as follows:
\begin{equation}
B'_{\textmd{FITA}}=\mathcal{A}(D'\circ_A (D_1\rightarrow B_1),\cdots,D'\circ_A (D_m\rightarrow B_m)),\end{equation}
where $\mathcal{A}$ is an $m$-ary aggregation function.

The other is First
Aggregate Then Infer (FATI). Concretely, the all fuzzy rules are aggregated into a single fuzzy relation, and then obtain the output
by composing an input $D'$ with the single fuzzy relation. With this scheme, the output can be expressed as follows:
\begin{equation}
B_{\textmd{FATI}}'=D'\circ_{A}\mathcal{A}(D_1\rightarrow B_1,\cdots, D_m\rightarrow B_m). \end{equation}

With the background as required in \cite{Zeng}, we assume that $D^i_{j}$ and $B_j$ are normal, continuous, complete and consistent pseudo-trapezoid-shaped which often form a Ruspini partition in the fuzzy rule base as a form (10). This means that the
fuzzy rules are complete and consistent.
\\{\bf Lemma 4.8} Assume that the number of  fuzzy rules
in (10) is greater than two. If the following conditions satisfy:

i. the operator $\rightarrow$
is chosen as a t-norm in inference algorithms (11) and (12),

ii. 0 is an annihilator of the aggregation functions $A$ and $\mathcal{A}$,

then $B'_{\textmd{FITA}}=B'_{\textmd{FATI}}\equiv0$.
\\{\bf Proof.} Since the fuzzy rules are complete and form a Ruspini partition, there exists
$j$ such that $D_j(\mathbf{x}_0)=0$ for an arbitrary given input $\mathbf{x}_0\in U^n$. Therefore, we have $B'_{\textmd{FITA}}(y)=\mathcal{A}(A(D'(\mathbf{x}_0), T(D_1(\mathbf{x}_0), B_1(y))),$  $\cdots,A(D'(\mathbf{x}_0),T(D_m(\mathbf{x}_0), B_m(y))))=\mathcal{A}(A(D'(\mathbf{x}_0), T(D_1(\mathbf{x}_0), $ $B_1(y))),\cdots, A(D'(\mathbf{x}_0),T(0,B_j(y))),\cdots,A(D'(\mathbf{x}_0),T(D_m(\mathbf{x}_0),B_m(y))))=0$.

Similarly, we can obtain $B'_{\textmd{FATI}}=0$.
\\{\bf Remark 5.} This result shows that we cannot choose aggregation functions having annihilator element 0 (Especially t-norms) to aggregate the inference results in Mamdani fuzzy system \cite{Mamdani}.

We can similarly obtain the following result.
\\{\bf Lemma 4.9} Assume that the number of  fuzzy rules
in (10) is greater than two. If the following conditions satisfy:

i. the operator $\rightarrow$
is chosen as a fuzzy implication in inference algorithms (11) and (12),

ii. 1 is an annihilator of the aggregation functions $A$ and $\mathcal{A}$,

then $B'_{\textmd{FITA}}=B'_{\textmd{FATI}}\equiv1$.
\\{\bf Remark 6.} This result shows that we cannot choose aggregation functions having annihilator element 1 (Especially t-conorms)  to aggregate the inference results in fuzzy logic controller \cite{Lee, Li2}.
\section{Similarity-based approximate reasoning with aggregation function}
\qquad In this section, we will extend the methods in \cite{Feng,li,Raha} using an aggregation function. And then investigate the validity of ASBR method.
Based on Definition 2.14, we say an inference scheme to solve the FMP and FMT problems satisfies

(GMP2$'$) $S(B',B)\leq S(B'',B)$ if $S(D',D)\leq S(D'',D)$.
\\{\bf Remark 7.} Obviously, (GMP2$'$) describes the fact that $B$ and $B'$ is more similar if $D'$ and $D$ is more similar.

Now, we extend the methods in \cite{Feng,li,Raha} to obtain the conclusion $B'$ of FMP problem. Inspired
by the ideas in \cite{Feng,Raha}, the following two modified conditional relations are considered to satisfy (AX1)-(AX3) proposed in \cite{Raha}:
$$R_1(D, B|D')(x, y)=A(S(D,D'), R(D,B)(x,y)),$$
$$ R_2(D, B|D')(x, y)=I(S(D,D'), R(D,B)(x,y)),$$
 where $A$ is an aggregation function and $I$ a fuzzy implication.

We further consider the following ways to construct the fuzzy relation $R(D, B)$:

Case 1. When $R(D, B)(x, y) =A(D(x), B(y))$, where $A$ is an aggregation function.

Case 2. When $R(D, B)(x, y) =I(D(x), B(y))$, where $I$ is a fuzzy implication.

Using sup-projection and inf-projection operations on $R_1(D, B|D')$ and $R_2(D, B|D')$ in case 1 and 2 respectively, we obtain
$$B'_1(y)=\sup_{x\in U}I(S(D,D'), A(D(x),B(y))),$$
$$B'_2(y)=\inf_{x\in U}I(S(D,D'), I(D(x),B(y))),$$
$$B'_3(y)=\sup_{x\in U}A(S(D,D'), A(D(x),B(y))),$$
$$B'_4(y)=\inf_{x\in U}A(S(D,D'), I(D(x),B(y))).$$
{\bf Lemma 5.1} Let $A$ be an aggregation function and $I$ a fuzzy implication. Suppose that $D$ is normal. For all corresponding inferred conclusions $B'_i(i=1,2,3,4)$, then we have

i. If $A$ has a left neutral element 1 and $I$ satisfies (NP), then $B\subseteq B'_1$;

ii. If $I$ is right-continuous and satisfies (NP), then $B\subseteq B'_2$;

iii. If $A$ has a left neutral element 0, then $B\subseteq B'_3$;

iv. If $A$ is right-continuous and has a left neutral element 0, then $B\subseteq B'_4$;
\\{\bf Proof.} i. Since $D$ is normal, there exists $x_0\in U$ such that $D(x_0)=1$. Thus, $B'_1(y)=\mathop{\sup}\limits_{x\in U}I(S(D,D'), A(D(x),B(y)))
\geq I(S(D,D'), A(D(x_0),B(y)))=\vspace{1mm}I(S(D,D'), A(1,B(y)))=I(S(D,D'),B(y))\geq I(1,B(y))=B(y)$.

ii-iv can be proved similarly to i.
\\{\bf Lemma 5.2} Let $I$ be a right-continuous fuzzy implication satisfying (NP). Suppose that $A$ has a left neutral element 1. If the inferred conclusion is determined by $B_1'$ or $B_2'$, then the method satisfies (GMP2$'$).
\\{\bf Proof.} Assume that $D', D''$ are two promises in FMP problem satisfying $S(D, D')\leq S(D, D'')$ and $B', B''$ are their corresponding conclusions. Then we have $B'_1(y)=\mathop{\sup}\limits_{x\in U}I(S(D,D'), A(D(x),$\vspace{1mm} $B(y)))
\geq \mathop{\sup}\limits_{x\in U}I(S(D,D''), A(D(x),B(y)))=B''_1(y)$. By Lemma 5.1, $B\subseteq B_1''\subseteq B_1'$ holds. This\vspace{1mm} implies $S(B, B_1')\leq S(B, B_1'')$.

We can similarly obtain $S(B, B_2')\leq S(B, B_2'')$ if $S(D, D')\leq S(D, D'')$.
\\{\bf Lemma 5.3} Let $A$ be an aggregation function and $I$  a fuzzy implication. If the conclusion of FMP problem is $B_3'$ or $B_4'$, then the method satisfies (GMP2).
\\{\bf Proof.} This result comes from the monotonicity of aggregation function.

For two fuzzy sets $D$ and $D'$, it is reasonable to assume that the measure of similarity is zero if and only if $\min(D(x),D'(x))=0$ holds for all $x\in U$. Therefore, we suppose $S(D,D^C)=0$ in order to discuss whether the above method satisfies (GMP3). And then we have the following results.
\\{\bf Lemma 5.4}  If the inferred conclusion is determined by $B_1'$ or $B_2'$, then the method satisfies (GMP3).
\\{\bf Proof.} Obviously.
\\{\bf Lemma 5.5} Let $A$ have a neutral element 0. Suppose that $D$ is normal. If the inferred conclusion is determined by $B_3'$, then the method satisfies (GMP3).
\\{\bf Proof.} Since 0 is a neutral element of $A$, $A(1,x)=1$ holds for any $x\in [0,1]$. For any $y\in [0,1]$, we have $B'_3(y)=\mathop{\sup}\limits_{x\in U}A(S(D,D), A(D(x),B(y)))=\mathop{\sup}\limits_{x\in U}A(1, A(D(x),$ $B(y)))=1$.\vspace{1mm}
\\{\bf Remark 8.} However, it is difficult to ensure $B_4'$ satisfying (GMP3).
\\{\bf Lemma 5.6} Let $A$ have a left neutral element 1 and $I$ satisfy (NP). Suppose that $D$ is normal. If the inferred conclusion is determined by $B_1'$, then the method satisfies (GMP4).
\\{\bf Proof.} Let $D'=D$.  By the monotonicity of $A$ and $I$, $I(S(D,D), A(D(x),B(y)))=I(1, A(D(x),$ $B(y)))=A(D(x),B(y))\leq A(1,B(y))=B(y)$ holds for any $x\in U$. This means $B'_1(y)=\mathop{\sup}\limits_{x\in U}I(S(D,D), A(D(x),B(y)))\leq B(y)$. According to Lemma 5.1, we have $B=B'_1$.\vspace{1mm}
\\{\bf Lemma 5.7} Let $I$ be right-continuous and satisfy (NP). Suppose that $D$ is normal. If the inferred conclusion is determined by $B_2'$, then the method satisfies (GMP4).
\\{\bf Proof.} This proof is similar to that of Lemma 5.6.
\\{\bf Lemma 5.8} Let $A$ have a left neutral element 1.  Suppose that $D$ is normal. If the inferred conclusion is determined by $B_3'$, then the method satisfies (GMP4).
\\{\bf Proof.} This proof is similar to that of Lemma 5.6.
\\{\bf Lemma 5.9} Let $A$ be right-continuous and have a left neutral element 1.  Suppose that $D$ is normal.
If the inferred conclusion is determined by $B_4'$, then the method satisfies (GMP4).
\\{\bf Proof.} This proof is similar to that of Lemma 5.6.
\section{Quintuple implications principle method of fuzzy inference with aggregation function}
\qquad This section will consider the AQIP method and its validity. Based on Lemma 2.15, we always suppose that the fuzzy implication is right-continuous with respect to the second variable in the rest of this paper.
\\{\bf Theorem 6.1} Let $I$ be a right-continuous with respect to the second variable fuzzy implication. If $I$ satisfies the condition $I(1,y)<1$ for all $y\in[0,1)$. Then the QIP solutions of FMP and FMT are as
follows:
\begin{equation}
B'(y)=\sup_{x\in U}A_I(A_I(D'(x),A_I(I(D'(x),D(x)),I(D(x),B(y)))),1),
\end{equation}
\begin{equation}D'(x)=\sup_{y\in V}A_I(A_I(D(x),A_I(I(D(x),B(y)),I(B(y),B'(y)))),1),\end{equation}
where $A_I$ is an aggregation function defined as Eq.(6).
\\{\bf Proof.} We only prove the QIP solution for FMP. The proof of FMT is similar. Since $I$ is right-continuous with respect to the second variable, the QIP solutions of FMP is unique by Lemma 2.15. Let $B'$ be defined as in Eq.(13). We firstly can verify that $B'$ can ensure that Eq.(1) takes its maximum 1. The right-continuous with respect to the second variable of $I$ implies that $A_I$ defined in Eq.(6) and $I$ satisfy (RP*) by Lemma 3.2. Then we have
$I(I(D(x),B(y)), I(I(D'(x),D(x)),
I(D'(x),B'(y))))=I(I(D(x),B(y)), I(I(D'(x),D(x)),I(D'(x),$ $\mathop{\sup}\limits_{x\in U}A_I(A_I(D'(x),A_I(I(D'(x),D(x)),I(D(x),B(y)))),1))))\geq1\Longleftrightarrow A_I(I(D(x),B(y)),1)$\vspace{1mm} $\leq I(I(D'(x),D(x)),I(D'(x),\mathop{\sup}\limits_{x\in U}A_I(A_I(D'(x),A_I(I(D'(x),D(x)),I(D(x),B(y)))),1)))\Longleftrightarrow A_I(I$\vspace{1mm} $(D'(x),D(x)),A_I(I(D(x),B(y)),1))\leq I(D'(x),\mathop{\sup}\limits_{x\in U}A_I(A_I(D'(x),A_I(I(D'(x),D(x)),I(D(x),$\vspace{1mm} $B(y)))),1)\Longleftrightarrow A_I(A_I(D'(x),A_I(I(D'(x),D(x)),I(D(x),B(y)))),1)\leq \mathop{\sup}\limits_{x\in U}A_I(A_I(D'(x),A_I(I($ \vspace{1mm}$D'(x),D(x)),I(D(x),B(y)))),1)$.

On the other hand, suppose that $C$ is an arbitrary fuzzy subset on $V$ such that
$I(I(D(x),B(y)),$ $I(I(D'(x),D(x)),
I(D'(x),C(y))))\equiv 1$ holds for any $x\in V$ and $y\in U$.  Since $A_I$ and $I$ satisfy (RP*), then
$$I(I(D(x),B(y)),I(I(D'(x),D(x)),
I(D'(x),C(y))))\equiv 1\qquad\ $$
$$\Longleftrightarrow A_I(I(D(x),B(y)),1)\leq I(I(D'(x),D(x)),
I(D'(x),C(y)))$$
$$\ \ \Longleftrightarrow A_I(I(D'(x),D(x)),A_I(I(D(x),B(y)),1))\leq
I(D'(x),C(y))$$
$$\quad\ \ \Longleftrightarrow A_I(D'(x),A_I(I(D'(x),D(x)),A_I(I(D(x),B(y)),1)))\leq C(y).$$
For $B'$ defined as in Eq.(13), this means that $B'(y)\leq C(y)$ holds for all $x\in U$ and $y\in V$.

The following example shows that QIP method for FMP does not satisfy (GMP1).
\\{\bf Example 6.2} Let $D=1/x_1+0.2/x_2+0.5/x_3$, $B=0.5/y_4+1/y_5$ and $D'=0.5/x_2+1/x_3+\vspace{2mm}$
$0.2/x_4$. If $I$ is chosen Goguen implication, that is, $I(x,y)=\left\{\begin{array}{cc}
          \frac{y}{x}\wedge 1 & x>0\\
          1& x=0
        \end{array}\right.$, then we have\vspace{1mm} $B'(y_5)=0.5<1=B(y_5)$.
\\{\bf Lemma 6.3} Let $I$ satisfy (NP) and (OP). For the solution of QIP method, then we have $B'\subseteq B$.
\\{\bf Proof.} We can assert that $A_I$ defined as Eq.(6) is commutative and associative similarly to Theorem 2.5.15 in \cite{Baczynski}.
Therefore, $A_I(D'(x),A_I(I(D'(x),D(x)),I(D(x),B(y))))\leq A_I(D(x),I(D(x),$ $B(y)))\leq B(y)$ holds. This implies that $B'\subseteq B$.
\\{\bf Lemma 6.4} Let $I$ satisfy (NP), (OP) and the equation $A_I(x,I(x,y))=x\wedge y$ for any $x,y\in [0,1]$.  Then  the solution of QIP method for FMP satisfies (GMP2).
\\{\bf Proof.} Obviously.

The following example shows that QIP method for FMP does not satisfy (GMP3).
\\{\bf Example 6.5} Let $D=1/x_1+0.2/x_2+0.5/x_3$ and $B=0.5/y_4+1/y_5$. We can compute $D^C=0.8/x_2+0.5/x_3+1/x_4+1/x_5$. If $I$ is chosen Goguen\vspace{1mm} implication, then we get $B'(y_1)=0<1$.
\\{\bf Lemma 6.6} Let $I$ satisfy (NP) and (IP).  If $D$ is normal, then the QIP method for FMP satisfies (GMP4).
\\{\bf Proof.} Let $D'=D$.  We have $I(D(x),D(x))=1$ for all $x\in U$. It is not difficult to see that $I(I(D(x),B(y)),I(I(D(x),D(x)),I(D(x),$ $B(y))))=1$ holds for all $x\in U$ and $y\in V$. According to quintuple implication principle for FMP, we have $B'(y)\leq B(y)$. Since $I$ satisfies (NP) and (IP), $A_I(1,x)=A_I(x,1)=x$ holds for all $x\in [0,1]$ by Theorem 3.4. Considering that $D$ is normal, there exists $x_0\in U$ such that $D(x_0)=1$. This means $B'(y)=\mathop{\sup}\limits_{x\in U}A_I(A_I(D(x),A_I(I(D(x),D(x)),$\vspace{1mm} $I(D(x),B(y)))),1)\geq A_I(A_I(1,A_I(I(1,1),I(1,B(y)))),1)$
$=B(y)$. Thus, $B'=B$.
\section{Discussion on three approximate reasoning methods}
\qquad In this section, we always assume that both $D'$ and $D^C$ are normal. For convenience, let ASBR$_i(i=1,2,3,4)$ denote the ASBR method to obtain $B'_i(i=1,2,3,4)$ in Section 5. From the above discussion, we can list together the three approximate reasoning methods and the GMP rules by which they satisfy as shown in Table 1. Notice that the satisfaction (dissatisfaction, respectively) of GMP rule is denoted by a $\surd$ ($\times$, respectively) in the column.

Clearly, it depends completely on the fuzzy implication and aggregation function whether the three approximate reasoning methods satisfy the GMP rules. Especially, the more properties of fuzzy implication and aggregation function are required in order to satisfy the GMP rules in the ASBR method and AQIP method. Therefore, it is not difficult to see that the ACRI methods should be a top priority of approximate reasoning according to the GMP rules. However, notice that the ASBR method and AQIP method  can effectively overcome the deficiency of ACRI method
\cite{Li1,Raha,Turksen,Zhou}. This implies that other properties (such as robustness, universal approximation capability etc.) should be utilized to measure the validity of aforementioned three approximate reasoning methods.
 $$\mbox{\bf{\small Table\ 1 \ three approximate reasoning methods and GMP rules}}$$
 \begin{center}
 \tabcolsep 0.05in
\begin{tabular}{cccccc}
 \toprule[1pt]
    Approximate reasoning methods & (GMP1) & (GMP2)&(GMP2$'$) &(GMP3)&(GMP4)\vspace{1mm}\\

   \midrule[0.75pt]
 ACRI method & $\surd$ & $\surd$&&$\surd$&$\surd$\vspace{1mm}\\
 ASBR$_1$ method & $\surd$ & &$\surd$&$\surd$&$\surd$\vspace{1mm}\\
 ASBR$_2$ method & $\surd$ & &$\surd$&$\surd$&$\surd$\vspace{1mm}\\
 ASBR$_3$ method & $\surd$ & $\surd$&&$\surd$&$\surd$\vspace{1mm}\\
 ASBR$_4$ method & $\surd$ & $\surd$&&$\times$&$\surd$\vspace{1mm}\\
 AQIP method & $\surd$ & $\times$&&$\times$&$\surd$\vspace{1mm}\\
   \bottomrule[1pt]
\end{tabular}
\end{center}\vspace{1mm}

Since the FMT is an extension of MT, as mentioned by Trillas in \cite{Trillas}, it is not trivial to further verify that whether the three approximate reasoning methods satisfy the following GMP rule:

(GMP5)\quad If $D'=B^C$, then $B'=D^C$.

Moreover, it is reasonable to involve some linguistic modifiers, such as very or little, in Premise 2 and conclusions of FMP and FMT problems. Therefore, we need to consider another GMP rule as follows.

(GMP6)\quad If $D'=m(D)$, then $B'=m(B)$, where $m$ is a modifier.
\section{Conclusions}
\qquad  Considering that the aggregation functions play a vital role in approximate reasoning and decision-making under imprecision or uncertainty, we firstly have utilized aggregation functions to construct three approximate reasoning methods. The validity of these three approximate reasoning  methods with aggregation functions has been further investigated. In our study, we have

(1)  Analyzed some properties of fuzzy implication
 generated by an aggregation function,

(2) Given the ACRI method with aggregation function,

(3) Studied the similarity-based approximate reasoning with aggregation function,

(4) Investigated the QIP solutions of FMP and FMT problems with aggregation function,

(5) Discussed the validity of three approximate reasoning methods aforementioned, respectively.

These results may act as a bridge between approximate reasoning and aggregation function. In the future, we
wish to investigate the capability of fuzzy inference system based on these methods. We also will apply them in prediction problems and decision making in real-life situation.
\section{Acknowledgement}  \qquad The authors would like to thank the anonymous referees and the Editor-in-Chief for their valuable comments.  This work was supported by the National Natural
Science Foundation of China (Grant No. 61673352).
\section{Compliance with ethical standards}
\vspace{2mm}
{\bf Conflict of interest} Author declares that he has no conflict of interest.
\vspace{1mm}
\\{\bf Human and animal rights} This article does not contain any studies
with human participants or animals performed by the authors.

\end{document}